\documentclass{amsproc}
\usepackage{amssymb}
\usepackage{amsfonts}
\usepackage{graphicx}
\newtheorem{theorem}{Theorem}

\newtheorem{examp}{Example}
\newtheorem{remar}{Remark}

\begin{document}
\newcommand*{\leftmapsto}{\mathbin{\reflectbox{$\longmapsto$}}}
\newcommand{\idi}{\mathfrak i}
\newcommand{\idig}{\underline{\mathfrak i}}
\newcommand{\idj}{\mathfrak j}
\newcommand{\idm}{\mathfrak m}
\newcommand{\idn}{\mathfrak n}
\newcommand{\un}{\underline{N}}
\newcommand{\unj}{\underline{N_1}}
\newcommand{\und}{\underline{N_2}}
\newcommand{\go}{\mathrm{germ}_0}
\newcommand{\wiw}{\mathrm{w}}
\newcommand{\ra}{\mathrm{rank}}
\newcommand{\orw}{\mathrm{ord}}
\newcommand{\au}{\mathrm{Aut}}
\newcommand{\re}{\mathrm{reg}}
\newcommand{\id}{\mathrm{id}}
\newcommand{\ima}{\mathrm{Im}}
\newcommand{\aop}{\er[t^1,\dots,t^k]}
\newcommand{\trk}{T^r_k}
\newcommand{\drk}{\mathbb D^r_k}
\newcommand{\dsl}{\mathbb D^s_l}
\newcommand{\ddk}{\mathbb D^2_k}
\newcommand{\ddd}{\mathbb D^2_2}
\newcommand{\ddj}{\mathbb D^2_1}
\newcommand{\djk}{\mathbb D^1_k}
\newcommand{\djl}{\mathbb D^1_l}
\newcommand{\dod}{\mathbb D \otimes \mathbb D}
\newcommand{\tdrk}{{\tilde{\mathbb D}}^r_k}
\newcommand{\tddk}{{\tilde{\mathbb D}}^2_k}
\newcommand{\tddd}{{\tilde{\mathbb D}}^2_2}
\newcommand{\bdrk}{{\bar{\mathbb D}}^r_k}
\newcommand{\bddk}{{\bar{\mathbb D}}^2_k}
\newcommand{\bddd}{{\bar{\mathbb D}}^2_2}
\newcommand{\ce}{C^{\infty}}
\newcommand{\nhk}{\er[s,t]/\langle st^2+s^4, s^2t+t^5\rangle+\idm^6}
\newcommand{\nhl}{\er[s,t]/\langle st^2+s^5, s^2t+t^5\rangle+\idm^6}
\newcommand{\mfm}{{\mathcal M\! f}_{\! m}}
\newcommand{\dex}{d \kern-.06em x}
\newcommand{\tut}{{\mathbf t}}
\newcommand{\oa}{\mathcal A}
\newcommand{\ek}{\mathcal E_k}
\newcommand{\coleq}{{\colon}{=}\;}
\newcommand{\idr}{\mathrm{id}_{\er}}
\newcommand{\osc}{\mathit{Osc}}
\newcommand{\er}{\mathbb R}
\newcommand{\de}{\mathrm d}
\newcommand{\ps}{\pi^s}
\newcommand{\psb}{\pi^s_b}
\newcommand{\aps}{{}_a^{}\!\pi^s_{}}
\newcommand{\apsb}{{}_a^{}\!\pi^s_b}
\newcommand{\po}{\pi^1}
\newcommand{\poo}{\pi^1_1}
\newcommand{\opo}{{}^{}_1\!\pi^1}
\newcommand{\pt}{\pi^2}
\newcommand{\pot}{\pi^1_2}
\newcommand{\porm}{\pi^1_{r-1}}
\newcommand{\tpo}{{}^{}_2\!\pi^1}
\newcommand{\iter}{\underbrace{T\dots T}_{\textrm{$r$--times}}}
\newcommand{\itrt}{\hbox{\space \raise 1.5mm\hbox{$\scriptstyle \circ \atop \textstyle T$}}\phantom{T}\!\!\!\!\!}
\newcommand{\itrtm}{\hbox{\space \raise 1.54mm\hbox{$\;{}_{{}_{\circ}} \atop \scriptstyle T$}}\phantom{T}\!\!\!\!\!}
\newcommand{\ka}{\mathcal K^A_{\er^m}}
\newcommand{\cg}{\mathcal G}
\newcommand{\cs}{{\bold{\mathcal S}}}
\newcommand{\ucs}{\underline{\cs}}
\newcommand{\fa}{\Phi_{\oa}}
\newcommand{\pa}{\Psi_{\oa}}
\newcommand{\fo}{\Omega^1 \er^m}
\newcommand{\dvj}{\frac{\partial}{\partial x^1}}
\newcommand{\dvm}{\frac{\partial}{\partial x^m}}
\newcommand{\rsa}{\rightsquigarrow}
\oddsidemargin 16.5mm
\evensidemargin 16.5mm

\thispagestyle{plain}

\begin{center}

{\Large\bf  LOOKING AT OSCULATING BUNDLES THROUGH THE SEMIHOLONOMITY EQUALIZATIONS
\rule{0mm}{6mm}\renewcommand{\thefootnote}{}%Enter at least one, but not more than 3 MSCs.
% First entered MSC will be a primary one, others (at most 2) will be secondary.
}

\vspace{1cc}
{\large\it Miroslav Kure\v{s}}

\vspace{1cc}
\parbox{24cc}{{\scriptsize
The paper glosses different forms of an introducing of higher order tangent-like functors,
especially functors derived from higher order nonholonomic tangent functors. 
A special attention is devoted to higher order osculating
bundles: their identification with higher order tangent bundles is demonstrated as the main result. Chiefly,
the paper is focused on the needful unification of concepts.
}}
\end{center}

\vspace{1cc}

\vspace{1.5cc}
\begin{center}
{\bf 1. MOTIONS AND VELOCITIES}
\end{center}

For clearness, we recall some basic concepts in a little bit of mechanics-like language.
Local differentiable maps between manifolds are differentiable maps defined on open subsets of a source manifold.
Further, our source manifold are always reals $\er$ (with a variable $t$) and our open subsets are (without loss of a generality) open real intervals.
Real intervals in question will be denoted by $I_1$, $I_2$, \dots.
If we take local maps from a real interval into a smooth manifold $M$, we talk about {\it local motions} on $M$.
Especially, we deal with local motions belonging to the same germ at $t_0\in\er$,
i.e. for $f\colon I_1\to M$, $t_0\in I_1$, $g\colon I_2\to M$, $t_0\in I_2$, the equality
$f(t_0)=g(t_0)=p\in M$ is satisfied.

Moreover, if $Y$ is a manifold fibered over $\er$, we can 
restrict, if it is needed, only to local {\sl sections} from an open real interval to $Y$.
In such a case, local motions are called {\it local time-transferring motions}.

Two local motions $f$ and $g$ belonging to the same germ
at $t_0\in\er$ can also belong to the same 1-jet ($r$-jet, respectively) at $u\in\er$. We write $j^1_{t_0}f=j^1_{t_0}g$ ($j^r_{t_0}f=j^r_{t_0}g$) and
the space of all 1-jets ($r$-jets) at $t_0\in\er$ is denoted by $J^1_{t_0}(\er,M)$ ($J^r_{t_0}(\er,M)$).
For $t_0=0\in\er$, we talk about the same {\it velocity} ({\it $r$-velocity}) in $p=f(0)=g(0)$
and $J^1_0(\er,M)=TM$ yields the {\it tangent bundle} ($J^r_0(\er,M)=T^rM$ yields the {\it $r$-th order tangent bundle}).
Nevertheless, we see that the choice $t_0=0$ (the choice of "true" zero) has only a formal character.

It is known that $TM$ is also a smooth manifold. Thus, {\sl first}, let us focus on velocities on the tangent bundle $TM$. Surely, we can again
take two local motions $F$, $G$ on $TM$:
the same velocity is considered in $F(0)=G(0)=(p,v)\in TM$. Forming the 1-jet, we see that we have no correspondence with the 
"installation" of $v$ now. 
Let us consider 
local coordinates $x^i\colon U\to\er^m$, $i=1,\dots,m$, $U\subseteq M$, $U\ni p$, $m=\dim M$. Then
we have the following expression in local coordinates: 
\begin{eqnarray*}
\er\supseteq I_1\ni 0 \stackrel{f}{\longmapsto} p\in M && \\
p \stackrel{(x^i)}{\longmapsto} \left(x^i(p)\right) &=& \left( x^i \left(f(0)\right)\right) \in\er^m \\
\er^m\ni\left(\frac{\partial x^i}{\partial t}(0) \right) &=& \left(y^i(v) \right) \stackrel{(y^i)}\leftmapsto v; \\
y^i \textrm{ are induced local coordinates}; &&\textrm{ all } (p,v) \textrm{ form } TM ;\\
\er\supseteq I_2\ni 0 \stackrel{F}{\longmapsto} (p,v)\in TM &&\\
(p,v) \stackrel{(x^i,y^i)}{\longmapsto} \left(\left(x^i,y^i\right)(p,v)\right)&=&
\left(\left(x^i,y^i\right)\left(F(0)\right)\right) \in\er^{2m} \\
\er^{2m}\ni\left(\left(\frac{\partial x^i}{\partial t},\frac{\partial y^i}{\partial t}\right)(0) \right)&=&\left(\left(X^i,Y^i\right)(P,V) \right) \stackrel{(X^i,Y^i)}\leftmapsto (P,V);\\
(X^i,Y^i) \textrm{ are induced local coordinates}; &&\textrm{ all } (p,v,P,V) \textrm{ form } TTM.\\
\end{eqnarray*}
So, $TTM=J^1_0\left(\er,J^1_0(\er,M)\right)$.

{\sl Second}, the {\it non-holonomic second order tangent bundle} is constructed by the following way. The space
of all 1-jets of local motions at all possible real points is denoted by $J^1(\er,M)$ 
(i.e. $J^1(\er,M)=\bigcup\limits_{t\in\er} J^1_t(\er,M)$).
As $J^1(\er,M)$ is a manifold fibered over $\er$, we can 
take local time-transferring motions now,
and, in particular, 1-jets of these time-transferring motions at $0$.
\begin{eqnarray*}
\er\supseteq I_1\ni t_0 \stackrel{f}{\longmapsto} p\in M &&\\
p \stackrel{(x^i)}{\longmapsto} \left(x^i(p)\right) &=&\left( x^i \left(f(t_0)\right)\right) \in\er^m \\
\er^m\ni\left(\frac{\partial x^i}{\partial t}(t_0) \right)&=&\left(y^i(v) \right) \stackrel{(y^i)}\leftmapsto v; \\
 y^i \textrm{ are induced local coordinates}; &&\textrm{ all } (t_0,p,v) \textrm{ form } J^1(\er,M);\\
\er\supseteq I_2\ni 0 \stackrel{\sigma}{\longmapsto} (0,p,v)\in J^1(\er,M) \\
(0,p,v) \!\stackrel{(\idr,x^i,y^i)}{\longmapsto}\! \left(\left(\idr,x^i,y^i\right)\!(0,p,v)\right)&=&
\left(\left(\idr,x^i,y^i\right)\!\left(\sigma(0)\right)\right) \!\in\!\{ 0 \}\!\!\!\times\!\!\er^{2m}\!\cong\!\er^{2m}\\
\er^{2m}\!\cong\!\{1\}\!\!\!\times\!\!\er^{2m}\!\ni\!\left(\left(\frac{\partial \idr}{\partial t},\frac{\partial x^i}{\partial t},\frac{\partial y^i}{\partial t}\right)(0) \right)&=&\left(\left(1,X^i,Y^i\right)(P,V) \right) \!\stackrel{(1,X^i,Y^i)}\leftmapsto \!(P,V); \\
(X^i,Y^i) \textrm{ are induced local coordinates}; &&\textrm{ all } (p,v,P,V) \textrm{ form } \tilde T^2M.\\
\end{eqnarray*}

So, we find that there is a slight difference between the second iterated tangent bundle and the non-holonomic second order tangent bundle.
Of course, this difference survives in general order.
However, iterated tangent functor and non-holonomic tangent functor are naturally equivalent.
The construction of the natural equivalence (in a general form) is noticed by Ivan Kol\'a\v{r} in \cite{KOL}; we refer
also to the paper \cite{KOV} of Ivan Kol\'a\v{r} and Raffaele Vitolo for more detailed description on pp.~4--5. 
We have mentioned only coordinate expressions grounds of the identification here.
In the paper, we follow routine to identify iterated tangent functor and non-holonomic tangent functor, however, we prefer
the iterated tangent functor in our notation.

\vspace{1.5cc}
\begin{center}
{\bf 2. PROJECTIONS IN ITERATED TANGENT BUNDLES}
\end{center}

From here, we will denote the iterated tangent functor $\iter$ by $\itrt^r$. This is a new notation. We believe that iterated tangent bundles deserves a special symbol. Intentionally, we will do the subsequent consideration for this bundle (and not for $\tilde T^rM$) as is perharps more frequent
and easily understandable.

Projections in $\itrt^r$ are well-known for number of decades, there were undoubtedly described already in 70's of last century, cf.
e.g. \cite{PAV} and references herein.
We introduce the following notation of projections in the iterated tangent
bundle $\itrt^r M$. (The notation was used by author in the paper \cite{KUN} and then, more precisely and with some basic properties of these
projections, in \cite{KUW}.) So, for every $s$, $0 < s\le r$, we denote by
$$
\ps\colon \itrt^s M\to M
$$
the canonical projection to the base.
Further, we denote
$$
\psb\coleq \ps_{\itrtm^b M} \colon \itrt^s\left(\itrt^b\right)\to \itrt^b M
$$
projection with $\itrt^s M$ as the base space,
$$
\aps\coleq \itrt^a \ps \colon \itrt^a \left(\itrt^ M \right) \to \itrt^a M
$$
induced projection originating by the posterior application of the
functor $\itrt^a$, and
$$
\apsb\coleq \itrt^a \ps_{\itrtm^b M}
$$ the general case
containing applications of both previous cases.
If $a$ or $b$ equal zero, we do not write them.

\begin{examp}\rm
For local coordinates $(x^i)$ on $M$, we obtain induced local coordinates $(x^i,y^i)$ on $\itrt^1 M$, $(x^i,y^i,X^i,Y^i)$ on $\itrt^2 M$ and
$(x^i,y^i,X^i,Y^i,\xi^i,\eta^i,\Xi^i,H^i)$ on $\itrt^3M$. We have
projections
\begin{eqnarray*}
\pi^3\colon \itrt^3 M\to M, \qquad && \pi^3(x^i,y^i,X^i,Y^i,\xi^i,\eta^i,\Xi^i,H^i)=(x^i),
\end{eqnarray*}
\begin{eqnarray*}
\pi^2_1\colon \itrt^3 M\to \itrt^1 M, \qquad && \pi^2_1(x^i,y^i,X^i,Y^i,\xi^i,\eta^i,\Xi^i,H^i)=(x^i,y^i),\\
{{}_1^{}\!\pi^2_{}}\colon \itrt^3 M\to \itrt^1 M, \qquad && {{}_1^{}\!\pi^2_{}}(x^i,y^i,X^i,Y^i,\xi^i,\eta^i,\Xi^i,H^i)=(x^i,\xi^i)
\end{eqnarray*}
and
\begin{eqnarray*}
\pi^1_2\colon \itrt^3 M\to \itrt^2 M, \qquad && \pi^2_1(x^i,y^i,X^i,Y^i,\xi^i,\eta^i,\Xi^i,H^i)=(x^i,y^i,X^i,Y^i),\\
{{}_1^{}\!\pi^1_1}\colon \itrt^3 M\to \itrt^2 M, \qquad && {{}_1^{}\!\pi^1_1}(x^i,y^i,X^i,Y^i,\xi^i,\eta^i,\Xi^i,H^i)=(x^i,y^i,\xi^i,\eta^i),\\
{{}_2^{}\!\pi^1_{}}\colon \itrt^3 M\to \itrt^2 M, \qquad && {{}_2^{}\!\pi^1_{}}(x^i,y^i,X^i,Y^i,\xi^i,\eta^i,\Xi^i,H^i)=(x^i,X^i,\xi^i,\Xi^i).
\end{eqnarray*}
We can obtain further projections by compositions, e.g. we have the projection $\poo\circ\tpo=\opo\circ\pot$ here,
which is not of the type $\apsb$. For details, see \cite{KUO}.
\end{examp}

\begin{remar}\rm
In the paper
of Ma\"{\i}do Rahula, Petr Va\v{s}\'{\i}k and Nicoleta Voicu \cite{RAV} is a completely different notation.
The projection $\pi_s$ there is our $\pi^1_{s-1}\colon \itrt^s M\to \itrt^{s-1} M$
and the projection $\rho_s$ there is our ${{}_{r-s}^{}\!\pi^1_{s-1}}\colon \itrt^r M\to \itrt^{r-1} M$.
Our notation is more general.
\end{remar}

\begin{remar}\rm
Elena Pavl\'{\i}kov\'{a} has used in \cite{PAV} still another notation. Her $j^s_r$ corresponds with our $\pi^{r-s}_s$ and
her ${}^p l^s_r$ corresponds with our ${}_p^{}\!\pi^{r-s}_{s-p}$. Projections ${}^p l^s_r$ are called (for $p\ge 1$) {\it lateral projections}
and some properties of these projections are derived in the cited paper.
\end{remar}

\vspace{1.5cc}
\begin{center}
{\bf 3. THE SEMIHOLONOMITY CONDITION AND THE OSCULATING BUNDLES}
\end{center}

Let $Z\in \itrt^r M$. We say that $Z$ is {\it prominent}, if the condition
$$
\pi^1_{r-1}(Z)=
{}^{}_{q-1}\!\pi^1_{r-q}(Z)
$$
for all $q$, $q=1,\dots,r$ is satisfied. (It is clear that a condition  
$\rho_1=\dots=\rho_r$ studied in \cite{RAV} is exactly the same.) This condition borrowed from the theory of nonholonomic jets is called the
{\it semiholonomity condition}.

The name "prominent" is only a working one. The prominent elements
of $\itrt^r M$ form a smooth manifold with fiber bundle structure over $M$. The obtained fibered manifold
is called (see \cite{RAV}) the {\it osculating bundle of manifold $M$} and denoted by
by $\osc^{r-1}M$. 

As we have seen before, the iterated tangent functor is different from the nonholonomic tangent functor, but there is a natural equivalence
between them. Using this equivalence, we can identify the bundle of prominent elements as the bundle of semiholonomic 
1-dimensional velocities of the order $r$. Nevertheless, semiholonomic 
1-dimensional velocities of the order $r$ are nothing but {\sl holonomic} velocities of the order $r$.
Hence the osculating bundle functor studied in \cite{RAV} is nothing but the higher order
tangent functor.

Thus, our result is the following.

\begin{theorem}
The osculating bundle functor $\osc^{r-1}$ 
works on the category of manifolds as objects and smooth maps as morphisms and it is naturally equivalent to the higher order tangent
func\-tor~$T^{r-1}$.
\end{theorem}
\textsl{To the completion of the proof.}
Almost everything concerning the proof of the theorem was already demonstrated above. We complete the proof by two comments.
First, every semiholonomic jet from 1-dimensional source manifold is automatically the holonomic one. Especially, this is evident from
local coordinate expressions: in general theory it was derived that nonholonomic jets which are identified after the equalizations
of projections differ only by displacements of zeros in subscripts, but no symmetrisation does not come yet (in general). Nevertheless, 
for a 1-dimensional source manifold an additional step of a symmetrisation is unreasoning because subscripts can have only one value. Second, the turn signal
for the observations in question can be a noticing dimensions of fibered manifolds: we recall that for $r\ge 2$ is $\dim T^{r-1}M = rm$
(with the fiber dimension $(r-1)m$) and cf. \cite{RAV} again.$\qquad\qquad\qquad\qquad\qquad\qquad\qquad\qquad\qquad\qquad\qquad\qquad\qquad\square$

\vspace{1.5cc}
\begin{center}
{\bf 4. FINAL REMARKS}
\end{center}

Thus, we hope it is clearer now, why, for instance, Wolfgang Bertram has remarked in his monographical work \cite{BER} (see Introduction, page 11) that 
{\sl "an osculating bundle of a vector bundle introduced by F. W. Pohl"} represents, {\sl "with great technical effort"}\dots
a construction of {\sl "a linear bundle which corresponds to $T^kF$"}.

Similarly, in projective geometry, the $k$-th osculating space is considered as the span of $\gamma(0)$, $\gamma^\prime(0)$, $\gamma^{\prime\prime}(0)$, \dots,
$\gamma^{(k)}(0)$ for a smooth parameterized curve $\gamma(t)$ (cf. the monograph of V. Ovsienko and S. Tabachnikov, \cite{OVT}), which suggests exactly
the same approach.

It follows that our result can be more or less known, maybe intuitively. We believe 
that this paper can be viewed as a contribution to the needful unification of different concepts.

\vspace{1.5cc}
\begin{center}
{\bf ACKNOWLEDGEMENTS}
\end{center}
The author was supported by GA \v{C}R, grant No. 201/09/0981.

%\section*{Bibliography}

%Fill the list of references in the alphabetical order here.
% You may use thebibliography  instead, and/or labels.
% Adopt the style like:

%\vspace{0.8cc}
%\newcounter{ref}
%\begin{list}{\small \arabic{ref}.}{\usecounter{ref} \leftmargin 4mm \itemsep -1mm}
%\item {\small {\sc J. A. Baker:} {\it Isometries in normed spaces.}
%Amer. Math. Monthly {\bf 78} (1971), 655--658.}
%\item {\small {\sc A. Marshall, I. Olkin}, "Inequalities: Theory of
%Majorization and Its Applications", Academic Press, New York, 1979.}
%\end{list}

\vspace{1cc}

 %Fill author(s) affiliation(s), address(es) and emails here:

{\small
\noindent

}
\begin{thebibliography}{9}


\bibitem{BER}
Bertram, W.,
{\it Differential Geometry, Lie Groups and Symmetric Spaces},
Memoirs of AMS, No. 900, 2008.

\bibitem{KOL}
Kol\'a\v{r}, I.,
{\it Bundle Functors of the Jet Type},
Differential Geometry and Applications, Proceedings of the 7th International Conference DGA 98, Brno, August 10-14, 1998, Masaryk University Brno, pp. 231--237 (1999).
      
\bibitem{KOV}
Kol\'a\v{r}, I., Vitolo, R.,
{\it Absolute Contact Differentiation on Submanifolds of Cartan Spaces},
Diff. Geom. Appl. {\bf 28}, No.1, 19--32 (2010).

\bibitem{KUN}
Kure\v{s}, M.,
{\it On the Simplicial Structure of Some Weil Bundles},
Supplem. ai Rendic. del Circ. Matem. di Pa\-ler\-mo~{\bf 63}, 131--140 (2000).

\bibitem{KUO}
Kure\v{s}, M.,
{\it On the Symmetrisation of Nonholonomic Jets},
Math. Proc. of Royal Irish Acad. {\bf 105}, No.A, 93--106 (2005).

\bibitem{KUW}
Kure\v{s}, M.,
{\it Weil Algebras of Generalized Higher Order Velocities Bundles},
Cont. Math. {\bf 288}, 358--362 (2001).

\bibitem{OVT}
Ovsienko, V., and Tabachnikov, S.,
{\it Projective Differential Geometry Old and New: from the Schwarzian Derivative to the Cohomology of Diffeomorphism Groups},
Cambridge University Press 2004.

\bibitem{PAV}
Pavl\'{\i}kov\'{a}, E., 
{\it Lateral Projections of Non-holonomic Jets},
Math. Slovaka {\bf 23}, No.2, 184--190 (1973).

\bibitem{RAV}
Rahula, M., Va\v{s}\'{\i}k, P., Voicu, N., 
{\it Tangent Structures: Sector-forms, Jets and Connections},
preprint.
\end{thebibliography}
\end{document}